\documentclass[11pt, oneside,onecolumn]{article}   	
\usepackage{geometry}                		
\usepackage{float}
\usepackage[caption = false]{subfig}
\usepackage{graphicx}				
							
\usepackage{amssymb, amsmath}
\usepackage{enumerate}

\newtheorem{thm}{Theorem}[section]

\newtheorem{prf}{Proof}[section]

\title{Ramanujan Sums as Derivatives}

\author{Devendra Kumar Yadav, Gajraj Kuldeep and S. D. Joshi}
\begin{document}

\maketitle

\title{\textbf {Abstract:-}} In 1918 S. Ramanujan defined a family of trigonometric sum now known as Ramanujan sums. In the last few years, Ramanujan sums have inspired the signal processing community. In this paper, we have defined an operator termed here as Ramanujan operator. In this paper it has been proved that these operator possesses properties of  first  derivative and second derivative with a particular shift. Generalised multiplicative property and new method of computing Ramanujan sums are also derived in terms of interpolation.\\

\textbf {Index terms}: Discrete Fourier Transform(DFT), Ramanujan sums(RS), Interpolation.

\section{Introduction} Ramanujan defined a trigonometric sum by taking $n$-th power of $q$-th primitive roots of unity ~\cite{bib:8}. A similar kind of sum can be seen in generation of DFT of a constant signal where constraint of primitive root is not used . \\
Until now RS are used for finding the period of the given signal. In this paper operator based on RS are defined which can give the detail portion of the signal. This makes it useful in find the edges of the given signal.

\section{Ramanujan Operator}

A new class of operators, termed here as Ramanujan operators, are defined using Ramanujan sums and it is also proved that these operators satisfies all the properties of a first derivative. 

Mathematically RS are defined as,
\begin{equation}
	c_q(n) = \sum_{\begin{subarray}{c}
	                         k=1\\
	                         (k, q) =1\\
	                         \end{subarray}}^{q}\exp^{\frac{j2\pi kn}{q}}
\end{equation}
where $(k,q)=1$ implies that k and q are relatively prime.

The Ramanujan sums $c_q(n)$ defined in equation $(1)$ are $q$ periodic. Few of them,for one period are given  below :
\begin{eqnarray*}
\begin{aligned}
	& c_1(n)= 1\\
	& c_2(n)= 1 , -1\\
	& c_3(n)= 2, -1 , -1\\
	& c_4(n)= 2, 0 , -2, 0\\
	\end{aligned}
\end{eqnarray*}

For a given $q$, $($q$ \ne 1)$ we define Ramanujan operator $\hat{R_q}$, as a linear, shift invariant operator whose kernel is defined as one period of $c_q(n)$. We denote the kernel of this operator as $\hat c_q(n)$. Mathematically,
\begin{eqnarray*}
\hat c_q(n) = \left(\begin{array}{cc}
											c_q(n) & 0 \le n \le q-1\\
											0 & otherwise
											\end{array}\right)
\end{eqnarray*}
Applying this operator on a signal $x(n)$,we get
\begin{eqnarray*}											
(\hat{R_q}x)(n) \triangleq \sum _{k}\hat c_q(k)x(n-k)
\end{eqnarray*}

Derivatives in digital domain are defined in terms of differences~\cite{bib:1}~\cite{bib:2}~\cite{bib:3}~\cite{bib:7}. In image processing first derivative is defined as a function that should possess the following properties.

\begin{enumerate}[1.]
	\item In areas of constant intensity derivative must be zero.
	\item The derivative must be nonzero at the onset of ramps or unit step.
	\item The derivative must be nonzero constant along ramps.
\end{enumerate}  
In addition to the properties $1$ and $2$ mentioned above, if the derivative is zero along ramps of constant slope then it is called as second derivative~\cite{bib:2}~\cite{bib:7}.
\begin{thm}
	For any given $q$ $($q$ \ne 1)$ Ramanujan operator $\hat{R_q}$ satisfies all the properties of first derivative(Hence can be used as a derivative operator).
\end{thm}
Some examples are $\hat{c_2}(n)=\{1,-1\}$,$\hat{c_3}(n)=\{2,-1,-1\}$.
\begin{prf}
	We will prove one by one that $\hat{R_q}$ satisfies all the three properties of first derivative.
	\begin{enumerate}[1.]
	\item For any $q$, the sum of $\hat c_q(n)$ for one period is zero. So in the areas of constant intensity it would trivially result in zero signal.
	\item Let us consider the unit step input to be of the form 
		
		\begin{eqnarray*}
			u(n-n_0) = 1 , n \ge n_0
		\end{eqnarray*}
		
		Applying $\hat{R_q}$ along the unit step we get,
		\begin{eqnarray*}
			f(n) = \sum_{l = 0}^{q-1}{u(n-n_0-l)\hat c_q(l)}
		\end{eqnarray*}
		
	For the case when $n$ is between $n_0$ and $n_0+q-1$ sum of $\hat c_q(n)$ will be nonzero, thereby giving a nonzero at the onset of the unit step input.
	
	\item Let us consider the ramp input to be of the form
	
		\begin{eqnarray*}
			r(n) = n 
		\end{eqnarray*}
		
		Applying $\hat{R_q}$ along the ramp we get,
		
		\begin{eqnarray*}
			f(n) = \sum_{l = 0}^{q-1}{(n-l)\hat c_q(l)}\\
			f(n) =\sum_{l = 0}^{q-1}{n\hat c_q(l)}-\sum_{l = 0}^{q-1}{l\hat c_q(l)}
		\end{eqnarray*}
		
		since $\hat c_q(l)$ is zero over one period. Therefore we get
		
		\begin{eqnarray*}
			f(n) =-\sum_{l = 0}^{q-1}{l\hat c_q(l)} = \sum_{\begin{subarray}{c}
	                         k=1\\
	                         (k, q) =1\\
	                         \end{subarray}}^{q}\frac{q}{(1-\exp^{\frac{j2\pi k}{q})}}
		\end{eqnarray*}
		
		This shows that $f(n)$ is constant which proves that the Ramanujan operator is nonzero constant along ramps.
	\end{enumerate}
	This concludes the proof of Theorem $2.1.$
\end{prf}
	
\begin{thm}
For a given odd $q$,$(q \ne 1)$ a new class of Ramanujan operators $\tilde{R_q}$ is defined which possess properties of second derivative. The kernel of these operators are represented by  $\tilde{c_q}(n)$ which is defined as first $q$ coefficients of $c_q(n)$ after giving a shift of $\frac{q-1}{2}$ to Ramanujan sum $c_q(n)$. 
Mathematically,
\begin{eqnarray*}
\tilde c_q(n) = \left(\begin{array}{cc}
											c_q(n-\frac{q-1}{2}) & 0 \le n \le q-1\\
											0 & otherwise
											\end{array}\right)
\end{eqnarray*}
\end{thm}
Some examples are $\tilde{c_3}(n)=\{-1,2,-1\}$,$\tilde{c_5}(n)=\{-1,-1,4,-1,-1\}$. 
\begin{prf}

From definition $c_q(n) = c_q(q-n)$. From definition $\tilde{c_q}(n)$ satisfies\\
\begin{eqnarray}
	\sum_{n = 0}^{q-1}{\tilde{c_q}(n)} = 0 \\
	\tilde{c_q}(i) = \tilde{c_q}(q-1-i)
\end{eqnarray}
It can be shown that this operator satisfies the first two properties of second derivative. Now we need to show that application of this operator on ramp with constant slope gives zero, for this operator to be a second derivative, 

	\begin{eqnarray*}
			f(n) = \sum_{l}{(n-l)\tilde{c_q}(l)}\\
			f(n) = \sum_{l = 0}^{q-1}{n\tilde{c_q}(l)}-\sum_{l = 0}^{q-1}{l\tilde{c_q}(l)}
	\end{eqnarray*}
	since $\tilde c_q(l)$ is zero over one period. Therefore we get
	\begin{eqnarray*}
			f(n) = -\sum_{l = 0}^{q-1}{l\tilde{c_q}(l)} = -(0\tilde{c_q}(0)+1\tilde{c_q}(1)+\hdots \\ +\frac{(q-1)}{2}\tilde{c_q}(\frac{q-1}{2})+\hdots +(q-1)\tilde{c_q}(q-1))
	\end{eqnarray*}	
	using equation $2$ and $3$, we get
	\begin{eqnarray*}
			f(n) = -(q-1)[\tilde{c_q}(0)+\tilde{c_q}(1)+\hdots +\frac{1}{2}\tilde{c_q}(\frac{q-1}{2})]=0
	\end{eqnarray*}
It can be further shown that $\sum_{l = 0}^{q-1}{(n-l)^2 \tilde{c_q}(l)}$ is a constant. This proves that Ramanujan operator $\tilde{R_q}$ acts as a second derivative.This conclude the proof of Theorem 2.2.
\end{prf}
In short, $c_q(n)$ acts as second derivative for a shift of $\frac{q-1}{2}$ when $q$ is odd. For other shifts it acts as a first derivative. When $q$ is even $c_q(n)$ acts as a first derivative for all shifts.

\section{Generalisation of multiplicative property and efficient computation of Ramanujan sums}
Ramanujan sums are multiplicative in nature.i.e.
\begin{eqnarray}
c_{pq}(n) = c_{p}(n)c_{q}(n)  \text{  if  } (p,q)=1
\end{eqnarray}
In this section multiplicative property of Ramanujan sums are generalised which means given any two shifted Ramanujan sums  their product also gives a shifted Ramanujan sum. A recursive way of computing Ramanujan sums is given in ~\cite{bib:8}. In this section we have proposed a non-recursive and efficient way of computing Ramanujan sums.
\begin{thm}
Generalisation of multiplication property: Given two arbitrary but fixed shifted Ramanujan sequences $c^{\alpha_1}_p(n)$ and $c^{\alpha_2}_q(n)$ where $\alpha_1$ and $\alpha_2$ are shift in the Ramanujan sequences , their product also forms a Ramanujan sequence.i.e.\\
\begin{eqnarray*}
	c^{\alpha_1}_p(n)c^{\alpha_2}_q(n) = \left\{\begin{array}{cc}
									c_{pq}(n-\alpha_1) & \text{ when } \alpha_1 = \alpha_2\\
									c_{pq}(n-\alpha_2p+\alpha_1q) & otherwise 
								\end{array}\right\}
\end{eqnarray*}
	
	where
	$c^{\alpha_1}_p(n) = c_p(n-\alpha_1)$, $c^{\alpha_2}_q(n) = c_q(n-\alpha_2)$,p$>$q and $(p,q)=1$ for $0 \le n \le pq-1$
\end{thm}

\begin{prf}
From definition we can write
\begin{eqnarray*}
	c^{\alpha_1}_p(n)c^{\alpha_2}_q(n) = \sum_{\begin{subarray}{c}
					                         k=1\\
	                         					(k, p) =1\\
	                         					\end{subarray}}^{p}\exp^{\frac{j2\pi k(n-\alpha_1)}{p}}\sum_{\begin{subarray}{c}
					                         l=1\\
	                         					(l, q) =1\\
	                         					\end{subarray}}^{q}\exp^{\frac{j2\pi l(n-\alpha_2)}{q}}
\end{eqnarray*}

	Let $W_q = exp^\frac{-j2\pi}{q}$, rewriting above equation,we get
	
\begin{eqnarray*}
	c^{\alpha_1}_p(n)c^{\alpha_2}_q(n) = \sum_{\begin{subarray}{c}
					                         k=1\\
	                         					(k, p) =1\\
	                         					\end{subarray}}^{p}W^{-k(n-\alpha_1)}_{p}
						\sum_{\begin{subarray}{c}
					                         l=1\\
	                         					(l, q) =1\\
	                         					\end{subarray}}^{q}W^{-l(n-\alpha_2)}_{q}
\end{eqnarray*}

\begin{eqnarray}
= \sum_{\begin{subarray}{c}
					                         k=1\\
	                         					(k, p) =1\\
	                         					\end{subarray}}^{p}
						\sum_{\begin{subarray}{c}
					                         l=1\\
	                         					(l, q) =1\\
	                         					\end{subarray}}^{q}W^{-n(kq+lp)}_{pq}W^{k\alpha_1}_{p}W^{l\alpha_2}_{q}
\end{eqnarray}

if $\alpha_1 = \alpha_2$ then from equation $5$ we get,
\begin{eqnarray*}
c^{\alpha_1}_p(n)c^{\alpha_1}_q(n) = c_{pq}(n-\alpha_1)
\end{eqnarray*}

Otherwise as we know that if $(p,q)=1$ then $(p-q,p)=1$. Therefore $(k(p-q),p)=1$. Rewriting equation $5$ using this we get,

\begin{eqnarray*}
	c^{\alpha_1}_p(n)c^{\alpha_2}_q(n) = \sum_{\begin{subarray}{c}
					                         k=1\\
	                         					(k, p) =1\\
	                         					\end{subarray}}^{p}
						\sum_{\begin{subarray}{c}
					                         l=1\\
	                         					(l, q) =1\\
	                         					\end{subarray}}^{q}W^{-n(kq+lp)}_{pq}W^{k(p-q)\alpha_1}_{p}W^{l(p-q)\alpha_2}_{q}
\end{eqnarray*}	

\begin{eqnarray*}
							= \sum_{\begin{subarray}{c}
					                         k=1\\
	                         					(k, p) =1\\
	                         					\end{subarray}}^{p}
						\sum_{\begin{subarray}{c}
					                         l=1\\
	                         					(l, q) =1\\
	                         					\end{subarray}}^{q}W^{-n(kq+lp)}_{pq}W^{(kqp-kq^2)\alpha_1}_{pq}W^{(lp^2-lqp)\alpha_2}_{pq}
\end{eqnarray*}	
					
	Replacing $W^{(kqp-kq^2)\alpha_1}_{pq} = W^{(-lqp-kq^2)\alpha_1}_{pq} = W^{-(lp+kq)\alpha_1q}_{pq}$ 

	Similarly $W^{(lp^2-lqp)\alpha_2}_{pq} = W^{(lp+kq)\alpha_2p}_{pq}$

\begin{eqnarray*}
	c^{\alpha_1}_p(n)c^{\alpha_2}_q(n) = \sum_{\begin{subarray}{c}
					                         k=1\\
	                         					(k, p) =1\\
	                         					\end{subarray}}^{p}
						\sum_{\begin{subarray}{c}
					                         l=1\\
	                         					(l, q) =1\\
	                         					\end{subarray}}^{q}W^{-n(kq+lp)}_{pq}W^{-(lp+kq)\alpha_1q}_{pq}\\W^{(lp+kq)\alpha_2p}_{pq}
\end{eqnarray*}	

\begin{eqnarray*}
	= \sum_{\begin{subarray}{c}
					                         k=1\\
	                         					(k, p) =1\\
	                         					\end{subarray}}^{p}
						\sum_{\begin{subarray}{c}
					                         l=1\\
	                         					(l, q) =1\\
	                         					\end{subarray}}^{q}W^{-(n+\alpha_1q-\alpha_2p)(kq+lp)}_{pq}
\end{eqnarray*}			
which is same as
$c^{\alpha_1}_p(n)c^{\alpha_2}_q(n) = c_{pq}(n-\alpha_2p+\alpha_1q)$

\end{prf}

\begin{thm}
For a prime $q$,
	$c_{q^l}(n) = q^{l-1}c_{q}(\frac{n}{q^{l-1}})$ where $l > 0$
\end{thm}

\begin{prf}
From definition of $c_q(n)$ \\
\begin{eqnarray*}
	c_{q^l}(n) = \sum_{\begin{subarray}{c}
	                         k=1\\
	                         (k, q) =1\\
	                         \end{subarray}}^{q^l}\exp^{\frac{j2\pi kn}{q^l}}
\end{eqnarray*}

Let $r_1,r_2,\hdots,r_s$ be the terms which are relatively prime to $q$. Therefore other terms which are relatively prime to $q^l$ are of the form

$q+r_1,q+r_2,\hdots,q+r_s$\\
$2q+r_1,2q+r_2,\hdots,2q+r_s$\\
\vdots\\
$(q^{l-1}-1)q+r_1,(q^{l-1}-1)q+r_2,\hdots,(q^{l-1}-1)q+r_s$.\\
Using this, rewriting the above expression of $c_{q^l}(n)$
\begin{eqnarray*}
c_{q^l}(n) = \exp^{\frac{j2\pi r_1n}{q^l}}(1+\exp^{\frac{j2\pi n}{q^{l-1}}}+\hdots+\exp^{\frac{j2\pi (q^{l-1}-1)n}{q^{l-1}}}) \hspace{0.5in} \\
	+\exp^{\frac{j2\pi r_2n}{q^l}}(1+\exp^{\frac{j2\pi n}{q^{l-1}}}+\hdots +\exp^{\frac{j2\pi (q^{l-1}-1)n}{q^{l-1}}}) \hspace{0.5in} \\
  \vdots 	\hspace{1.5in}\\
	+\exp^{\frac{j2\pi r_sn}{q^l}}(1+\exp^{\frac{j2\pi n}{q^{l-1}}}+\hdots +\exp^{\frac{j2\pi (q^{l-1}-1)n}{q^{l-1}}}) \hspace{0.5in}\\
\end{eqnarray*}	

Since
$(1+\exp^{\frac{j2\pi n}{q^{l-1}}}+\exp^{\frac{j2\pi 2n}{q^{l-1}}}+\hdots +\exp^{\frac{j2\pi (q^{l-1}-1)n}{q^{l-1}}}) = q^{l-1}$ when $q^{l-1}|n$ and $0$ otherwise.\\
Therefore
\begin{eqnarray*}
	c_{q^l}(n) = (\exp^{\frac{j2\pi r_1n}{q^l}}+\exp^{\frac{j2\pi r_2n}{q^l}}+\hdots+\exp^{\frac{j2\pi r_sn}{q^l}})q^{l-1}\\
	=q^{l-1}\sum_{\begin{subarray}{c}
	                         k=1\\
	                         (k, q) =1\\
	                         \end{subarray}}^{q}\exp^{\frac{j2\pi kn}{q^l}}
	= q^{l-1}c_{q}(\frac {n}{q^{l-1}})
\end{eqnarray*}	

\end{prf}

\begin{thm}
For arbitrary positive integer N, Ramanujan sum can be expressed in terms of Ramanujan sums of its prime factors. 
\end{thm}

\begin{prf}
Let N be any arbitrary number which can be represented as $N$  = $p_1^{r_1}p_2^{r_2} \hdots p_m^{r_m}$
Using multiplicative property $c_N(n)$ can be written as
\begin{eqnarray*}
	c_N(n) = c_{p_1^{r_1}}(n)c_{p_2^{r_2}}(n) \hdots c_{p_m^{r_m}}(n)
\end{eqnarray*}
Using theorem $3.2$ $c_N(n)$ can be rewritten as
\begin{eqnarray*}
	c_N(n) = p_1^{r_1-1}p_2^{r_2-1} \hdots p_m^{r_m-1} c_{p_1}(\frac {n}{p_1^{r_1-1}})c_{p_2}(\frac {n}{p_2^{r_2-1}}) \hdots c_{p_m}(\frac {n}{p_m^{r_m-1}})
\end{eqnarray*}
This proves theorem $3.3$.
\end{prf}
It shows that RS of higher order can be seen as interpolated version of RS of lower order.
\section{Concluding Remarks}
In this paper, Ramanujan operator has been defined.  It has been proved that these operators possesses properties of  first  derivative and second derivative also if given a particular shift.  Apart from these operators generalised multiplicative property and new method of computing Ramanujan Sums are also derived. This can be used for image processing which will be explored in future.

\section*{Acknowledgment}

The authors would like to thank Dr. Dhananjoy Dey for supportive guidance.

\end{document}